\documentclass[a4paper,12pt,leqno]{article}
\usepackage{amsmath,amsfonts,amsthm,amssymb,dsfont}
\usepackage{graphicx}
\usepackage[alphabetic]{amsrefs}
\usepackage[OT4]{fontenc}
\usepackage{enumerate}

\long\def\symbolfootnote[#1]#2{\begingroup%
\def\thefootnote{\fnsymbol{footnote}}\footnote[#1]{#2}\endgroup}

\newtheorem{thm}{Theorem}[section]

\newtheorem{prop}[thm]{Proposition}
\newtheorem{cor}[thm]{Corollary}

\theoremstyle{definition}
\newtheorem{rem}[thm]{Remark}
\newtheorem{defin}[thm]{Definition}
\newtheorem{constr}[thm]{Construction}
\newtheorem{ex}[thm]{Example}
\newcommand{\VH}{\mathcal{VH}}
\newcommand{\Flat}{\mathbb{E}^2_\triangle}

\begin{document}

\begin{center}
\large\bfseries Square complexes and simplicial nonpositive curvature
\end{center}

\begin{center}\bf
Tomasz Elsner$^a$\symbolfootnote[1]{Partially supported by MNiSW grant
N N201 541 738.} and
Piotr Przytycki$^b$\symbolfootnote[2]{Partially supported by MNiSW grant
N N201 541 738 and the Foundation for Polish Science.}
\end{center}

\begin{center}\it
$^a$ Mathematical Institute, University of Wroc\l aw,\\
Plac Grunwaldzki 2/4, 50-384 Wroc\l aw, Poland\\
\emph{e-mail:}\texttt{elsner@math.uni.wroc.pl}
\end{center}

\begin{center}\it
$^b$ Institute of Mathematics, Polish Academy of Sciences,\\
 \'Sniadeckich 8, 00-956 Warsaw, Poland\\
\emph{e-mail:}\texttt{pprzytyc@mimuw.edu.pl}
\end{center}

\begin{abstract}
\noindent
We prove that each nonpositively curved square $\VH$-complex can be turned functorially into a locally $6$-large simplicial complex of the same homotopy type. It follows that any group acting geometrically on a CAT(0) square $\VH$-complex is systolic. In particular the product of two finitely generated free groups is systolic, which answers a question of Daniel Wise. On the other hand, we exhibit an example of a compact non-$\VH$ nonpositively curved square complex, whose fundamental group is neither systolic, nor even virtually systolic.
\end{abstract}



\section{Introduction}

In this note we compare nonpositively curved square $\VH$-complexes (introduced in~\cite{WPhD}) and locally $6$-large simplical complexes (introduced in \cite{JS}). First we describe locally $6$-large and systolic complexes. The definitions we use are taken from \cite{JS2}, with a slight modification allowing simplices in a locally $6$-large simplicial complex not to be embedded. Nevertheless, the definition of a systolic complex coincides with the one in \cite{JS2}.

\begin{defin}
A \emph{generalised simplicial complex} is a set $\mathcal{S}$ of affine simplices together with a set $\mathcal{E}$ (closed under compositions) of affine embeddings of simplices of $\mathcal{S}$ onto the faces of simplices of $\mathcal{S}$ (attaching maps), such that for any proper face $\tau$ of any simplex $\sigma\in\mathcal{S}$ there is \emph{precisely one} attaching map onto $\tau$.

A \emph{(generalised) simplicial map} between generalised simplicial complexes is a set of affine maps commuting with the attaching maps and mapping each source simplex onto one of the target simplices.

The \emph{geometric realisation} of a generalised simplicial complex $(\mathcal{S},\mathcal{E})$ is the quotient space $\mathcal{S}/\mathcal{E}$. The quotient map of a generalised simplicial map is the \emph{geometric realisation} of such a map. We will abuse the language and not distinguish between simplicial complexes or simplicial maps and their geometric realisations.

The \emph{link} of a simplex $\sigma$ in a complex $X=(\mathcal{S},\mathcal{E})$ is the (generalised) simplicial complex $X_\sigma=(\mathcal{S}_\sigma,\mathcal{E}_\sigma)$ where the set $\mathcal{S}_\sigma$ is obtained by taking for each attaching map $\phi_{\sigma,\tau}\colon\sigma\to\tau$ the maximal subsimplex of $\tau$ disjoint from the image of $\sigma$
and $\mathcal{E}_\sigma$ is the set of restrictions of the maps in $\mathcal{E}$.
\end{defin}

Subsequently, we refer to a generalised simplicial complex simply as a \emph{simplicial complex} and use the phrase \emph{simple simplicial complex} when referring to a standard simplicial complex (in which simplices are embedded and the intersection of two simplices, if non-empty, is a single simplex).

\begin{defin}
A simplicial complex is \emph{simple} if it does not contain an edge joining a vertex to itself, or a pair of simplices with the same boundary (e.g. a double edge). A simple simplicial complex is $\emph{flag}$ if any complete subgraph (a clique) of its $1$-skeleton spans a simplex.

A \emph{cycle without diagonals} in a simplicial complex $X$ is an embedded simplicial loop such that there are no edges in $X$ connecting a pair of its nonconsecutive vertices.

A simplicial complex is \emph{locally $6$-large} if all of its vertex links are flag and do not contain cycles of length $4$ or $5$ without diagonals.
A connected and simply connected locally $6$-large simplicial complex is called \emph{systolic} (i.e. systolic complexes are the universal coverings of connected locally $6$-large complexes).

A group admitting a geometric action on a systolic complex is called \emph{systolic}.
\end{defin}

The original definition of local 6-largeness in \cite{JS2} requires that we check the flagness and the absence of short cycles without diagonals for the link at \emph{any} simplex. However, for higher-dimensional simplices it is a direct consequence of those properties for the links at the vertices.

Similarly as for simplicial complexes, we allow cells in square complexes not to be embedded. The formal definition of a (generalised) square complex is the same as of a generalised simplicial complex, except for putting \emph{vertices, edges and squares} in place of \emph{simplices}. The only thing that needs to be rephrased is the definition of the link.

\begin{defin}
The \emph{link} at a vertex $v$ of a (generalised) square complex $X=(\mathcal{S},\mathcal{E})$ is a 1-dimensional (generalised) simplicial complex $X_v=(\mathcal{S}_v,\mathcal{E}_v)$ (a graph), where $\mathcal{S}_v$ is obtained by taking for each attaching map $\phi_{v,\sigma}\colon v\to\sigma$ the vertex of $\sigma$ opposite to $v$ (if $\sigma$ is an edge) or the diagonal of $\sigma$ opposite to $v$ (if $\sigma$ is a square) and $\mathcal{E}_v$ is the set of restrictions of the maps in $\mathcal{E}$.

A square complex is called a \emph{$\VH$-complex} if its $1$-cells can be partitioned into two classes $V$ and $H$ called \emph{vertical} and \emph{horizontal} edges, respectively, and the attaching map of each square alternates between the edges of $V$ and $H$. In other words, the link at each vertex is a bipartite graph with independent sets of vertices coming from edges of $V$ and $H$.
\end{defin}

Note that the link of a $\VH$-complex at a vertex may have double edges.

\begin{defin}
A square complex is \emph{nonpositively curved} (or \emph{locally CAT(0)}) if the link at any vertex does not contain embedded combinatorial cycle of length less than $4$. For a $\VH$ complex this reduces to the property that there are no double edges in the links at vertices.
\end{defin}

For a general definition of CAT(0) and nonpositively curved spaces (not needed in our article) see \cite{BH}. Note only, that a simply connected space which is nonpositively curved is CAT(0) (\cite[Theorem 4.1]{BH}).

\begin{ex}
The product of two trees is a CAT(0) $\VH$-complex. If a group acts freely by isometries on the product of two trees and preserves the coordinates, then the quotient square complex is a nonpositively curved $\VH$-complex.
\end{ex}

The paper is divided into two parts. In Section \ref{sectiontwo} we provide a functorial construction turning a nonpositively curved square $\VH$-complex into a locally 6-large simplicial complex of the same homotopy type (in particular turning a CAT(0) $\VH$-complex into a systolic complex). The main application of the construction is:

\begin{thm}[see Corollary \ref{fundgp}]\label{intro}
The fundamental group of a compact nonpositively curved $\VH$-complex is systolic.
\end{thm}

The first application of Theorem \ref{intro} is the answer to a question posed by Daniel Wise in \cite{W}:

\begin{cor}
The product of two finitely generated free groups is systolic.
\end{cor}

We also obtain a series of consequences of Theorem \ref{intro} by applying it to the examples of nonpositively curved $\VH$-complexes (some with exotic properties) given by Daniel Wise in \cite{WPhD}.

\begin{cor}[{compare \cite[Corollary 2.8]{WPhD}}]
The fundamental group of an alternating knot complement is systolic.
\end{cor}

\begin{cor}[{compare \cite[Theorem 5.5]{WPhD}}]
There exists a systolic group, which is not residually finite.
\end{cor}

One can arrange for even a stronger property:

\begin{cor}[{compare \cite[Theorem 5.13]{WPhD}}]
There exists a systolic group, which has no finite-index subgroups.
\end{cor}

In Section 3 we show that the $\VH$-hypothesis in Theorem \ref{intro} is necessary:

\begin{thm}[see Theorem \ref{nonvirt}]\label{intro-nonVH}
There exists a compact non-$\VH$ nonpositively curved square complex, whose fundamental group is not systolic, nor even virtually systolic.
\end{thm}

\textbf{Acknowledgements.}
We thank Daniel Wise for motivating us, for his suggestions and discussions.

\section{Nonpositively curved $\VH$-complexes are systolic}\label{sectiontwo}

Our main construction yields a way of turning a nonpositively curved $\VH$-complex into a locally $6$-large simplicial complex.

\begin{constr}
\label{constr}

Let $X$ be a $\VH$ complex with the sets $E_V$ and $E_H$ of vertical and horizontal edges, respectively. Denote by $V$ and $S$ the sets of vertices and squares of $X$, respectively. We construct an associated simplicial complex $X^*$ called the \emph{simplexification} of $X$, which has the same homotopy type as $X$.

First we divide each vertical edge $e\in E_V$ in two and subdivide each square $s\in S$ into six triangles, as in the Figure \ref{simplex}(a), obtaining a (generalised) simplicial complex $\hat X$ (a triangulation of $X$). The vertices of $\hat X$ (which will correspond to the vertices of $X^*$) are in bijective correspondence with the elements of $V\cup E_V\cup S$. We denote those vertices by $v^*$, $e^*$, $s^*$, for $v\in V$, $e\in E_V$, $s\in S$, respectively.

\begin{figure}[ht!]
\centering
\includegraphics{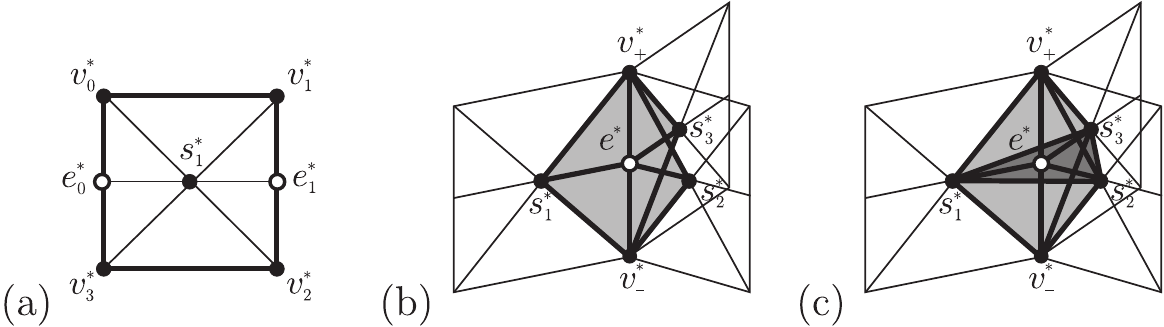}
\caption{(a) the triangulation of $X$ \hskip.5cm (b) $\hat Y_e\subset \hat X$ \hskip.5cm (c) $Y^*_e\subset X^*$}
\label{simplex}
\end{figure}

The link of $\hat X$ at a vertex $e^*$ is isomorphic to the suspension of a set of $n$ points, where $n$ is the number of squares $s\in S$ with a vertical edge $e$ (counted with multiplicities, i.e. a square with both vertical edges equal to $e$ is counted twice). The union $\hat Y_e$ of all the simplices of $\hat X$ containing the vertex $e^*$ is isomorphic to the suspension of an $n$-pod, where some pairs of vertices may be identified.

The complex $X^*$ is obtained from $\hat X$ by attaching simplices $\sigma^+_v=v^*_+e^*s^*_1\dots s^*_n$ and $\sigma^-_v=v^*_-e^*s^*_1\dots s^*_n$ for each vertex $e^*$, where $v_+$ and $v_-$ are the endpoints of the edge $e\in E_V$ and $s_1,\dots,s_n$ are the squares adjacent to the vertical edge $e\in E_V$ (counted with multiplicities).

The link of $e^*$ in $X^*$ is the suspension of an $(n-1)$-simplex and the union $Y^*_e$ of all the simplices of $X^*$ containing the vertex $e^*$ is isomorphic to the suspension of an $n$-simplex, where some pairs of vertices may be identified.
\end{constr}

\begin{prop}\label{homotopy}
A $\VH$ square complex $X$ and its simplexification $X^*$ have the same homotopy type.
\end{prop}

\begin{proof}
As the triangulation $\hat X$ of $X$ (defined in Construction \ref{constr}) embeds into $X^*$, we only need to prove that $X^*$ deformation retracts onto $\hat X$. Since for distinct $e_0,e_1\in E_V$ we have $Y^*_{e_0}\cap Y^*_{e_1}\subset \hat X$ it is enough to show that for any $e\in E_V$ the complex $Y^*_e$ deformation retracts onto $Y^*_e\cap \hat X=\hat Y_e$.

If $Y^*_e$ is a simple complex (i.e. the suspension of the simplex with vertices $s^*_1,\dots,s^*_n$,$e^*$), then denoting by $S$ the suspension and by $C$ the cone operator, we have
	$$\hat Y_e=S(C(\{s^*_1,\dots,s^*_n\}))\subset S(C(\sigma(s^*_1,\dots,s^*_n)))=Y^*_e,$$
where $\sigma(s^*_1,\dots,s^*_n)$ is the simplex with vertices $s^*_1,\dots,s^*_n$.

Consider the retraction $r:C(\sigma(s^*_1,\dots,s^*_n))\to C(\{s^*_1,\dots,s^*_n\})$ defined to be the affine extension of the map from the first barycentric subdivision, which preserves the subcomplex $C(\{s^*_1,\dots,s^*_n\})$ and maps the barycentres of the remaining simplices to the cone vertex. It is easy to see that $r$
can be extended to a deformation retraction. By suspending the deformation retraction, we obtain a deformation retraction from $Y^*_e$ onto $\hat Y_e$.

If $Y^*_e$ is not simple, then it is a quotient space of the suspension of a simplex obtained by identifying some pairs of vertices. In that case the deformation retraction from $Y^*_e$ onto $\hat Y_e$ is the quotient of the map described above.
\end{proof}

\begin{rem}
\label{functorial}
Note that Construction~\ref{constr} is functorial. Namely, let $f\colon X\rightarrow Y$ be a \emph{combinatorial map} between $\VH$ complexes (i.e.\  mapping cells onto cells of the same dimension, in our case mapping edges to edges and squares to squares). Assume also that $f$ preserves the sets of vertical and horizontal edges. Then $f$ induces a canonical combinatorial map $f^*\colon X^*\rightarrow Y^*$. Moreover, we have $(f\circ g)^*=f^*\circ g^*$ and $id^*=id$.
In particular, if $f$ is invertible (in other words is a \emph{combinatorial isomorphism}; it induces an isometry between the geometric realisations), then so is $f^*$. Finally, note that if a group $G$ acts properly (cocompactly, geometrically) on $X$, then its induced action on $X^*$ is also proper (cocompact, geometric).
\end{rem}

We are now ready for our main result.

\begin{thm}
\label{main}
If $X$ is a nonpositively curved $\VH$ complex, then its simplexification $X^*$ is locally $6$-large.
\end{thm}

Before giving the proof, we list a few consequences, obtained by applying Proposition~\ref{homotopy} and Remark~\ref{functorial}.

\begin{cor}
\label{first cor}
If $X$ is a CAT(0) $\VH$ complex, then its simplexification $X^*$ is systolic. If $G$ acts geometrically on $X$, then $G$ is systolic.
\end{cor}

There are two notable applications of Corollary~\ref{first cor}.

\begin{cor}[Theorem \ref{intro}]\label{fundgp}
The fundamental group of a compact nonpositively curved $\VH$ complex is systolic.
\end{cor}

The second application promotes Wise's aperiodic flat construction (\cite[Construction 7.1]{WPhD}) into the systolic setting.

\begin{defin}\label{sysflat}
A \emph{flat} in a systolic complex $X$ is a subcomplex $\Flat\subset X$ which is isomorphic to the equilaterally triangulated plane (the triangulation with 6 triangles adjacent to each vertex) and whose 1-skeleton is isometrically embedded into $X^{(1)}$ (with the combinatorial metric).
\end{defin}

\begin{cor}
There exists compact a locally 6-large simplicial complex, whose universal cover (which is systolic) contains a flat, which is not the limit of a sequence of periodic flats.
\end{cor}

It remains to prove our main result.

\begin{proof}[Proof of Theorem~\ref{main}]

We need to check that the link of $X^*$ at any vertex is flag and does not contain cycles of length 4 or 5 without diagonals. It is immediate for any vertex $e^*$, where $e\in E_V$, as the link of $e^*$ is the suspension of a simplex.

The link of $X$ at a vertex $s^*$, $s\in S$ has the form of two suspensions of simplices, $S\sigma_{m-1}$ and $S\sigma_{n-1}$ ($m$ and $n$ being the numbers of squares adjacent to the vertical edges of $s$), whose top and bottom vertices are connected by an edge (the case $m=n=3$ is depicted in Figure \ref{link}(b)). It is clear that it is flag and any cycle without diagonals in that link has length at least 6.

\begin{figure}[ht!]
    \centering
    \includegraphics{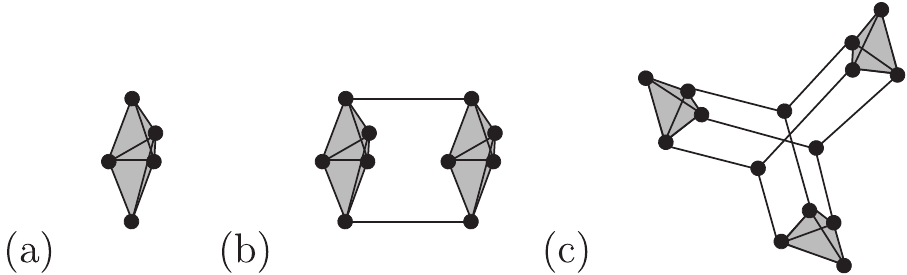}
    \caption{Sample link of $X^*$ at a vertex \ (a) $e^*$ \ (b) $s^*$ \ (c) $v^*$}\label{link}
\end{figure}

Now let $L$ be the link of $X^*$ at $v^*$, $v\in V$.
Then $L$ is the union of:
\begin{itemize}
\item a set of simplices (one $n_e$-simplex for each vertical edge $e\in E_V$ with an endpoint $v$, where $n_e$ is the number of squares $s\in S$ adjacent to $e$, counted with multiplicities) and
\item a set of $m_e$-pods (one $m_e$-pod for each horizontal edge $e\in E_H$ issuing from $v$, where $m_e$ is the number of squares adjacent to $e$, counted with multiplicities),
\end{itemize}
where each endpoint of each $m_e$-pod is identified with a different vertex of one of the simplices (Figure \ref{link}(c) shows the link $L$ in the case when the neighbourhood of $v\in V$ is the product of two tripods). It is clear that $L$ is flag, and any cycle without diagonals in $L$ needs to pass through at least two $m_e$-pods and two simplices, which makes its length at least 6.
\end{proof}

\section{Examples of nonpositively curved square complexes which are not systolic}

In the next part of the paper we show that our theorem cannot be improved to include all nonpositively curved square complexes. Namely, we construct an example of a compact nonpositively curved square complex, whose fundamental group is not systolic. Later, we use that example to show a compact nonpositively curved square complex, whose fundamental group is neither systolic, nor even virtually systolic (Theorem \ref{intro-nonVH}).

\begin{figure}[ht!]
    \centering
    \includegraphics{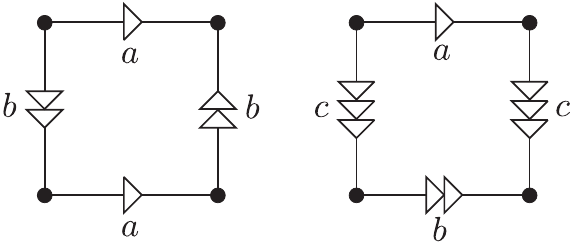}
    \caption{The nonpositively curved square complex $K$}\label{K}
\end{figure}

\begin{figure}[ht!]
    \centering
    \includegraphics{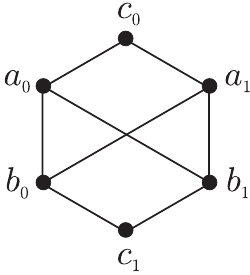}
    \caption{The link of $K$ at the only vertex}\label{Klink}
\end{figure}

Let $K$ be the square complex presented in Figure \ref{K}, built up of two squares. It has only one vertex and the link at this vertex is shown in Figure \ref{Klink}. Thus we see that $K$ is a nonpositively curved square complex, but not a $\VH$-complex. We will show that $\pi_1(K)$ is not a systolic group.

\begin{thm}\label{nonVH}
The group $\pi_1(K)$ is not systolic.
\end{thm}

\begin{proof}
The group $$\pi_1(K)=\left<a,b,c\;|\;ba=ab^{-1},\;a=cbc^{-1}\right>$$ is an HNN-extension of the fundamental group of a Klein bottle, so it has a subgroup $H=\left<a,b\right>$, which is isomorphic to the fundamental group of a Klein bottle, in particular is virtually $\mathbb{Z}^2$.

Suppose $\pi_1(K)$ is systolic, i.e. acts geometrically on some systolic simplicial complex $X$. As a corollary from the systolic flat torus theorem (precisely by Corollary 6.2(1) together with Theorem 5.4 in \cite{E}) we have that $H$, as a virtually $\mathbb{Z}^2$ group, acts properly on a systolic flat in $X$ (see Definition \ref{sysflat}). If the fundamental group of a Klein bottle $\left<a,b\;|\;ba=ab^{-1}\right>$ acts properly (by combinatorial isomorphisms) on an equilaterally triangulated plane $\Flat$ (shown in Figure \ref{flat}), then the axis of the glide reflection $a$ is $l$ and the direction of the translation $b$ is $k$ or vice versa.

\begin{figure}[ht!]
    \centering
    \includegraphics[scale=1.25]{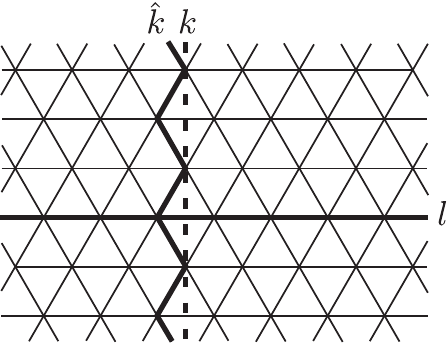}
    \caption{$\Flat$}\label{flat}
\end{figure}

The elements $a^2$ and $b^2$ act by translations on $\Flat$ (with axes $k$ and $l$). The 1-skeleton of $\Flat$ with the combinatorial metric is isometrically embedded into the 1-skeleton of $X$ (by Definition \ref{sysflat}), so the lines $\hat k$ and $l$ (marked in Figure \ref{flat}) are invariant geodesics (in the 1-skeleton) for $a^2$ and $b^2$. By \cite[Proposition 3.10]{E2} the geodesic $l$ is quasi-convex in the 1-skeleton of $X$ equipped with the combinatorial metric (i.e. any geodesic in $X^{(1)}$ with both endpoints on $l$ is contained in the $\delta$-neighbourhood of $l$, for some universal $\delta$). The geodesic $\hat k$ is clearly not quasi-convex (every point of $\Flat$ lies on some geodesic with both endpoints on $\hat k$).

Since $a^2=cb^2c^{-1}$, the translation $a^2$ has two invariant geodesics: $\hat k$ and $c(l)$ (or $l$ and $c(\hat k)$). Two invariant geodesics of an isometry acting by a translation on both of them are at finite Hausdorff distance, so either both $\hat k$ and $c(l)$ are quasi-convex, or none of them is (\cite[Proposition 3.11]{E2}). That contradicts the fact that $l$ is quasi-convex, while $\hat k$ is not.
\end{proof}

As we have just shown, the fundamental group of $K$ is not systolic, however it is virtually systolic (there is a 2-leaf covering $\tilde K$, which is a $\VH$-complex, so $\pi_1(\tilde K)$ is systolic by Theorem \ref{main}). Now we use the complex $K$ to construct a square complex $S$, whose fundamental group is not even virtually systolic.

Let $E$ be the compact nonpositively curved $\VH$-complex which has no connected finite coverings, constructed by Wise in \cite[Theorem 5.13]{WPhD}. Let $\sigma$ be any loop in $E$ consisting entirely of horizontal edges. We can subdivide the complex $K$ such that all loops $a$, $b$ and $c$ have the same combinatorial length as the loop $\sigma$. Now we define $\bar E$ and $\bar{\bar E}$ to be two copies of $E$ and let $$S=(E\cup \bar E\cup \bar{\bar E})\cup K/\sim,$$ where $\sim$ is the identification of $\sigma$, $\bar\sigma$ and $\bar{\bar\sigma}$ with $a$, $b$ and $c$, respectively. Then $S$ is a nonpositively curved (non-$\VH$) square complex.

\begin{thm}[Theorem \ref{intro-nonVH}]\label{nonvirt}
The group $\pi_1(S)$ is not virtually systolic, where $S$ is the nonpositively curved square complex defined above.
\end{thm}

\begin{proof}
We first argue that $\pi_1(S)$ is not systolic. Since
$$\pi_1(S)=\pi_1(K)*_{a=\sigma}\pi_1(E)*_{b=\bar\sigma}\pi_1(\bar E)*_{c=\bar{\bar\sigma}}\pi_1(\bar{\bar E})$$
is an amalgam product, the inclusion $K\subset S$ induces injection $\pi_1(K)\to\pi_1(S)$. To conclude that $\pi_1(S)$ is not systolic, we can recall the fact that a finitely presented subgroup of a torsion-free systolic group is systolic itself (\cite{W}), while $\pi_1(K)$ is not systolic (Theorem \ref{nonVH}). An equivalent way of arriving to that conclusion is to repeat for $S$ the argument used for $K$ in the proof of Theorem \ref{nonVH}.

To prove that $\pi_1(S)$ is not virtually systolic, we show that it has no finite-index subgroups (i.e. $S$ has no connected non-trivial finite coverings). Let $p:\tilde S\to S$ be a connected finite covering. Since $E\subset S$ has no connected non-trivial finite coverings, $p^{-1}(E)$ is a disjoint union of copies of $E$. In particular, any lift $\tilde a$ of the loop $a$ has the same length as $a$. The same holds for the loops $b$ and $c$. As $a$, $b$ and $c$ together with the three copies of $\pi_1(E)$ generate $\pi_1(S)$, that implies that $p$ is a trivial covering.
\end{proof}

\end{document}